\documentclass[12pt]{amsart}
\usepackage{amssymb}
\usepackage{amsxtra}
\usepackage{stmaryrd}
\usepackage{isolatin1}    
\usepackage[all]{xy}
\usepackage{graphics}

\newtheorem{thm}{Theorem}[section]
\newtheorem{lem}[thm]{Lemma}
\newtheorem{cor}[thm]{Corollary}
\newtheorem{prop}[thm]{Proposition}
\newtheorem{ex}[thm]{Example}

\newtheorem*{prob*}{Open problem}

\theoremstyle{definition}

\newtheorem{defi}[thm]{Definition}

\theoremstyle{remark}

\newtheorem{rem}[thm]{Remark}
\newtheorem*{rem*}{Remark}

\DeclareMathOperator{\Hom}{Hom}

\newcommand{\kringel}{\mathbin{\raise1pt\hbox{$\scriptstyle\circ$}}}
\newcommand{\pkt}{\mathbin{\raise0pt\hbox{$\scriptstyle\bullet$}}}

\newcommand{\C}{\mathbb{C}}

\newcommand{\N}{\mathbb{N}}

\newcommand{\R}{\mathbb{R}}

\newcommand{\Z}{\mathbb{Z}}

\newcommand{\End}{\mathop{\rm End}}

\newcommand{\Der}{\mathop{\rm Der}}

\newcommand{\Aut}{\mathop{\rm Aut}}

\newcommand{\Lie}{\mathop{\rm Lie}}

\newcommand{\La}{\mathfrak{a}}

\newcommand{\Lg}{\mathfrak{g}}
\newcommand{\Lh}{\mathfrak{h}}

\newcommand{\Ll}{\mathfrak{l}}
\newcommand{\Ln}{\mathfrak{n}}
\newcommand{\Lm}{\mathfrak{m}}

\newcommand{\Lp}{\mathfrak{p}}

\newcommand{\Lr}{\mathfrak{r}}
\newcommand{\Ls}{\mathfrak{s}}

\newcommand{\Lu}{\mathfrak{u}}
\newcommand{\Lv}{\mathfrak{v}}

\newcommand{\CC}{\mathcal{C}}

\newcommand{\CF}{\mathcal{F}}
\newcommand{\CG}{\mathcal{G}}

\newcommand{\CL}{\mathcal{L}}

\newcommand{\CN}{\mathcal{N}}
\newcommand{\CO}{\mathcal{O}}
\newcommand{\CR}{\mathcal{R}}

\newcommand{\im}{\mathop{\rm im}}
\newcommand{\abs}[1]{\lvert#1\rvert}

\newcommand{\al}{\alpha}
\newcommand{\be}{\beta}

\newcommand{\ep}{\varepsilon}

\newcommand{\la}{\lambda}

\newcommand{\ov}{\overline}

\newcommand{\ra}{\rightarrow}

\renewcommand{\phi}{\varphi}

\begin{document}

\title[]{Contractions of Lie algebras and algebraic groups}
\author[D. Burde]{Dietrich Burde}

\address{Fakult\"at f\"ur Mathematik\\
Universit\"at Wien\\
  Nordbergstrasse 15\\
  1090 Wien}

\date{\today}
\email{dietrich.burde@univie.ac.at}


\begin{abstract}
Degenerations, contractions and deformations of various algebraic structures
play an important role in mathematics and physics.
There are many different definitions and special cases of these
notions. We try to give a general definition which unifies these notions and shows the 
connections among them. Here we focus on contractions of Lie algebras and
algebraic groups.
\end{abstract}

\maketitle

\section{Contractions, degenerations and deformations of Lie algebras}

\subsection{Basic definitions and properties}

The notion of Lie algebra and Lie group contractions was
first introduced by I.E. Segal \cite{SEG} and E. In\"on\"u, E.P. Wigner \cite{IW}.
The usual definition of a continuous contraction of a Lie algebra is as follows.

\begin{defi}
Let $V$ be a vector space over $\R$ or $\C$ and $g\colon (0,1]
\ra GL(V)$ be a continuous function. Let $[,]$ be a Lie bracket
on $V$. A parametrized family of  Lie brackets on $V$ is defined by
\[
[x,y]_{\ep}=g_{\ep}([g_{\ep}^{-1}(x), g_{\ep}^{-1}(y)]).
\]
If the limit
\[
\llbracket x,y\rrbracket =\lim_{\ep \to 0} [x,y]_{\ep}
\]
exists, then $\llbracket , \rrbracket$ is a Lie bracket on $V$ and
$(V,\llbracket , \rrbracket)$ is called a {\it contraction} of $(V,[,])$.
\end{defi}

For $0<\ep\le 1$ the Lie algebras $(V,[,]_{\ep})$ are all isomorphic
to $(V,[,])$. Hence to obtain a {\it new} Lie algebra via contraction one needs
$\det (g_{\ep})=0$ for $\ep=0$. This is a necessary condition, but not a 
sufficient one. \\[0.2cm]
A contraction can be viewed as a special case of a so called degeneration.
Let $V$ be an $n$-dimensional vector space over a field $k$. Denote by
$\CL_n(k)$ the {\it variety of Lie algebra laws}. This is the set of all possible
Lie brackets $\mu$ on $V$. 
$\CL_n(k)$ is an algebraic subset of the affine variety
$\Lambda^2 V^{\ast} \otimes V$ of all alternating bilinear maps
from $V\times V$ to $V$. 
For a fixed basis $(x_1,\ldots,x_n)$ of $V$ a Lie bracket $\mu$ is
determined by the point $(c_{ij}^r)\in k^{n^3}$ of structure constants with
\begin{equation*}
\mu(x_i,x_j) = \sum_{r=1}^n c_{ij}^r x_r
\end{equation*}
satisfying the polynomial conditions
\begin{align*}
0 & =c_{ij}^r + c_{ji}^r, \\
0 & =\sum_{r=1}^n (c_{ij}^r c_{lr}^s+c_{jk}^r c_{ir}^s+c_{ki}^r c_{jr}^s).
\end{align*}
for $1\le i<j<k\le n,\; 1\le s\le n$, given by skew-symmetry and Jacobi's
identity. 
The general linear group $GL_n(k)$ acts on $V$, and hence on $\CL_n(k)$ by:
\begin{equation*}
(g\cdot\mu)(x,y)=g(\mu(g^{-1}x, g^{-1}y))
\end{equation*}
for $g\in GL_n(k)$ and $x,y\in V$.  
Denote by $O(\mu)$ the orbit of $\mu$ under this action, and
by $\ov{O(\mu)}$ the closure of the orbit with respect to the Zariski
topology. The orbits in $\CL_n(k)$ correspond to isomorphism classes of
$n$-dimensional Lie algebras.

\begin{defi}
Let $\la,\mu \in \CL_n(k)$ be two Lie algebra laws. We say that $\la$ degenerates
to $\mu$, if $\mu \in \ov{O(\la)}$. This is denoted by $\la \ra_{\rm deg} \mu$.
\end{defi}

A degeneration is called {\it trivial} if $\la \cong \mu$, that is, if
$\mu \in O(\la)$. 

\begin{rem}
Any irreducible component $\CC$ of $\CL_n(k)$ containing $\mu$ also contains
all degenerations of $\mu$. Indeed, we have  $O(\mu)\subset \CC$ so that
$\ov{O(\mu)}$ is contained in $\CC$, since  $\CC$ is closed. \\
A Lie algebra law $\mu$ is {\it rigid}, if its orbit $O(\mu)$ is open
in $\CL_n(k)$. Then $\ov{O(\mu)}$ is an irreducible component of $\CL_n(k)$.
There are only finitely many irreducible components in each dimension.
\end{rem}

\begin{rem}\label{1.3}
If $G$ is an algebraic group and $X$ is an algebraic variety over an 
algebraically closed field $K$, with regular action, then any orbit $G(x)$, 
$x\in X$ is a smooth algebraic variety, open in its closure $\ov{G(x)}$. Its 
boundary $\ov{G(x)}\setminus G(x)$ is a union of orbits of strictly lower
dimension. Each orbit $G(x)$ is a constructible set, hence $\ov{G(x)}$ coincides
with the closure $\ov{G(x)}^d$ in the standard topology. This can be found
in Borel's book \cite{BOR}, see the closed orbit lemma.
\end{rem}

Denote by $K$ an algebraically closed extension of $k$.

\begin{prop}
Degeneration in $\CL_n(K)$ defines a partial order on the orbit space of $n$-dimensional
Lie algebra laws by $O(\mu)\le O(\la) \iff \mu \in \ov{O(\la)}$.
\end{prop}

\begin{proof}
The relation is clearly reflexive. The transitivity follows from the fact that
$O(\la)\subseteq \ov{O(\mu)} \iff \ov{O(\la)}\subseteq \ov{O(\mu)}$.
Finally, antisymmetry follows from the fact, that any orbit in this case
is open in its closure, see remark $\ref{1.3}$. 
\end{proof}

The above order relation on the orbit space is represented by the so called
Hasse diagram. We repeat that degeneration is transitive: $\la \ra_{\rm deg} \mu$ and 
$\mu \ra_{\rm deg} \nu$ imply that $\la \ra_{\rm deg} \nu$.

\begin{lem}
Each Lie algebra contraction is a Lie algebra degeneration.
\end{lem}

\begin{proof}
Suppose $\la$ contracts to $\mu$. The subset $g_t \cdot\la$ of $O(\la)$ is 
parametrized by $t$. Therefore every polynomial
function vanishing on $O(\la)$ also vanishes on the Lie algebra laws in $t$, where
$t$ is replaced by $0$, hence on $\mu$. 
Therefore $\mu$ belongs to the Zariski closure of $O(\la)$.
\end{proof}

We may view $g_t$ formally as an element in $GL_n(k(t))$, where $k(t)$
is the field of fractions of the polynomial ring $k[t]$.

\begin{ex}
Every law $\la \in \CL_n(k)$ contracts to the abelian law 
$\la_0\in \CL_n(k)$.
\end{ex}

We have $\la_0(x,y)=0$, and with $g_t=t^{-1}I_n$ we have
$\la \ra_{\rm deg} \la_0$ since
$$(g_t \cdot \la)(x,y)=t^{-1}\la(tx,ty)=t\la(x,y).$$
Indeed, the limit of $g_t \cdot \la$ for $t\to 0$ equals $\la_0$. 
Some algebras like $\Lh_3\oplus k^{n-3}$, where $\Lh_3$ is the $3$-dimensional Heisenberg 
Lie algebra, can only degenerate to the abelian Lie algebra of the same dimension,
see \cite{LA}.
Given two Lie algebra laws $\la,\mu \in \CL_n(k)$ it is sometimes quite
difficult to see whether there exists a degeneration $\la \ra_{\rm deg} \mu$.
It is helpful to obtain some necessary conditions for the existence of a
degeneration.
In some sense one can say that $\la \ra_{\rm deg} \mu$ implies that
$\mu$ is ``more abelian'' than $\la$. A much finer condition is that
the dimensions of the cohomology spaces cannot decrease.

\begin{prop}
Let $\la \ra_{\rm deg} \mu$ a non-trivial degeneration over $\C$. Then we have
for all $i\in \N_0$:\\
\begin{align*}
\dim O(\la) & > \dim O(\mu)\\
\dim \Der \la  & < \dim \Der \mu\\
\dim \la^i & \ge \dim \mu^i \\
\dim \la^{(i)} & \ge \dim \mu^{(i)} \\
\al(\la) & \le \al (\mu) \\
{\rm rank} (\la) & \le {\rm rank} (\mu) \\
\dim Z(\la) & \le \dim Z(\mu) \\
\dim H^i(\la) & \le \dim H^i(\mu) \\
\dim H^i(\la,\la) & \le \dim H^i(\mu,\mu)\\
\end{align*}
where $\al(\la)$ denotes the maximal dimension of an abelian subalgebra
of $\la$, and
\begin{align*}
\la^0 & =\la^{(0)}=\la, \\
\la^i & =[\la,\la^{i-1}], \\
\la^{(i)}& =[\la^{(i-1)},\la^{(i-1)}].
\end{align*}
\end{prop}
For a proof see \cite{NEP} and the references given there. 
The first claim follows from Borel's closed orbit lemma, see remark $\ref{1.3}$.
Note that $O(\mu)$ can be identified with $GL(V)/\Aut (\mu)$, so
that
\begin{align*}
\dim O(\mu) & = \dim GL(V)-\dim \Aut (\mu) \\
 & = n^2-\dim \Der (\mu).
\end{align*}
This shows the second claim. The other claims rely also on the following lemma.

\begin{lem}
Let $G$ be a reductive algebraic group over $\C$ with Borel subgroup $B$. If
$G$ acts regularly on an affine variety $X$, then for all $x\in X$,
$\ov{G\cdot x}=G\cdot (\ov{B\cdot x})$.
\end{lem}

One can also show that the above proposition is also valid for $k=\R$. But then
additional arguments are needed. \\[0.2cm]
It is already quite interesting to investigate the varieties
$\CL_n(k)$ and the orbit closures over the complex numbers
in small dimensions.
\begin{ex}
For $n=2$ we have 
$$\CL_2(\C)=\ov{O(\Lr_2(\C))}=O(\Lr_2(\C)) \cup O(\C^2)$$
where $\Lr_2(\C)$ is the non-abelian algebra. 
\end{ex}

The only non-trivial degeneration is given by $\Lr_2(\C)\ra_{\rm deg} \C^2$.
The orbit of $\Lr_2(\C)$ is open. There is no Lie algebra law
degenerating to $\Lr_2(\C)$ in $\CL_2(\C)$.\\

\begin{ex}
The variety $\CL_3(\C)$ is the union of two irreducible components
$\CC_1$ and $\CC_2$.
\end{ex}

The component $\CC_1$ consists of the Lie algebras of trace zero,
i.e., where the linear form  ${\rm \mbox{tr ad}}(x)$ vanishes:
$$\CC_1=\ov{O(\Ls\Ll_2(\C))}=O(\Ls\Ll_2(\C)) \cup O(\Lr_{3,-1}(\C))
\cup O(\Ln_3(\C)) \cup O(\C^3)$$
The component $\CC_2$ consists of the solvable Lie algebras:
$$\CC_2=\CR_3(\C)=\cup_{\al}O(\Lr_{3,\al}(\C)) \cup O(\Lr_3(\C))\cup
O(\Lr_2(\C)\oplus \C)\cup O(\Ln_3(\C) \cup O(\C^3)$$
We have $\CC_1\cap \CC_2=\ov{O(\Lr_{3,-1}(\C))}$ and $\dim \CC_1=\dim \CC_2=6$.\\[0.5cm]
The classification of all orbits and their orbit closures in
$\CL_3(\C)$ is given as follows:
\vspace*{0.5cm}
\begin{center}
\begin{tabular}{c|c|c}
 $\Lg$ & Lie brackets & $\ov{O(\Lg)}$  \\
\hline
$\C^3$ & $-$ & $\C^3$ \\
$\Ln_3(\C)$ & $[e_1,e_2]=e_3$ & $\Ln_3(\C),\;\C^3$ \\
$\Lr_2(\C) \oplus  \C$ & $[e_1,e_2]=e_2$ & $\Lr_2(\C)\oplus\C,\;\Ln_3(\C),\;\C^3$ \\
$\Lr_3(\C)$ & $[e_1,e_2]=e_2,\, [e_1,e_3]=e_2+ e_3$ & $\Lr_3(\C),\; \Lr_{3,1}(\C),\; \Ln_3(\C),\;\C^3$ \\
$\Lr_{3,\al}(\C)$ & $[e_1,e_2]=e_2, \,[e_1,e_3]=\al e_3,\, \al\in I$
& $\Lr_{3,\al}(\C),\;\Ln_3(\C),\;\C^3 $\\
$\Lr_{3,-1}(\C)$ & $[e_1,e_2]=e_2, \,[e_1,e_3]=-e_3$
& $\Lr_{3,-1}(\C),\;\Ln_3(\C),\;\C^3 $\\
$\Lr_{3,1}(\C)$ & $[e_1,e_2]=e_2, \,[e_1,e_3]=e_3$ & $\Lr_{3,1}(\C),\;\C^3$ \\
$\Ls \Ll_2 (\C)$ & $[e_1,e_2]=e_3, [e_1,e_3]=-2 e_1,[e_2,e_3]=2 e_2$ &
$\Ls \Ll_2 (\C),\;\Lr_{3,-1}(\C),\;\Ln_3(\C),\;\C^3$
\end{tabular}
\end{center}

\vspace*{0.5cm}
Here for $\al,\be \neq 0$ we have $\Lr_{3,\al}(\C)\cong \Lr_{3,\be}(\C)$ if and only if
$\al=\be$ or $\be=\al^{-1}$. Let $I$ denote the set of $\al\in \C$ satisfying
$0<\abs{\al}\le 1$, and, if $\abs{\al}=1$, then $\al=e^{i\theta}$ with 
$\theta\in [0,\pi]$. The following Hasse diagram shows all essential degenerations (that is,
all the other degenerations are combinations of these) in
$\CL_3(\C)$, see \cite{BU10}:

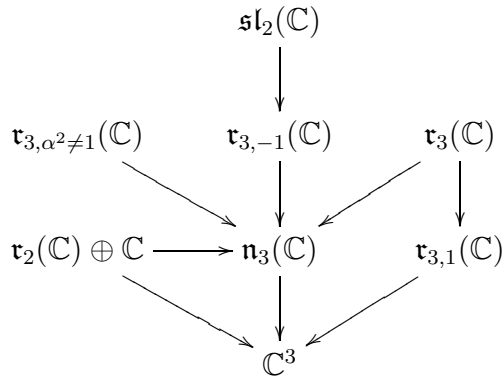
$$
\begin{xy}
\xymatrix{
 &  \Ls \Ll_2(\C) \ar[d] &  \\
\Lr_{3,\alpha^2 \ne 1}(\C)\ar[rd]  &  \Lr_{3,-1}(\C)\ar[d]  & \Lr_{3}(\C)\ar[ld]\ar[d] \\
\Lr_{2}(\C)\oplus \C \ar[rd] \ar[r]  & \Ln_3(\C) \ar[d] &  \Lr_{3,1}(\C)\ar[ld] \\
 & \C^3 &
}
\end{xy}
$$

For $n=4$ the classification of orbit closures is already quite complicated.
For details see \cite{BU10}, \cite{AGA}, \cite{CA1}. We have the following result:

\begin{prop}
The variety $\CL_4(\C)$ is the union of $4$ irreducible components
$\CC_i$, $i=1,\ldots,4$ as follows:
\begin{align*}
\CC_1 & = \ov{O(\Ls\Ll_2(\C)\oplus \C)} \\
\CC_2 & = \ov{O(\Lr_2(\C)\oplus \Lr_2(\C))}\\
\CC_3 & = \ov{\cup_{\al,\be}O(\Lg_4(\al,\be))}\\
\CC_4 & = \ov{\cup_{\al}O(\Lg_5(\al))}
\end{align*}
\end{prop}
Here $\Lg_4(\al,\be)$ has Lie brackets
\begin{align*}
[e_1,e_2]& = e_2, \\
[e_1,e_3]& = e_2+\al e_3,\,[e_1,e_4]=e_3+\be e_4,\\
\end{align*}
and $\Lg_5(\al)$ has Lie brackets
\begin{align*}
[e_1,e_2]& =e_2, \,[e_1,e_3]=e_2+\al e_3,\\
[e_1,e_4]& =(\al +1)e_4,\,[e_2,e_3]=e_4.
\end{align*}

The components are of dimension $12$, i.e., $\dim \CC_i=12$.
The number of open orbits equals $2$; indeed, the Lie algebras
$\Ls\Ll_2(\C)\oplus \C$ and $\Lr_2(\C)\oplus \Lr_2(\C)$
are rigid.\\
For computations of orbit closures for nilpotent Lie algebras (of dimension
$n\le 7$) see \cite{BU9}, \cite{BU18}. \\

\begin{defi}
Let $(\Lg,[\> , \, ])$ be a Lie algebra over $k$ and
$g,h\in \Lg,\;\phi_k \in \Hom (\Lambda ^2 \Lg,\Lg)$.
A {\it formal deformation} of $\Lg$ over $k[[t]]$ is a power series
\begin{equation*}
[g,h]_t:=[g,h]+\sum_{k\ge 1}\phi_k(g,h)t^k,
\end{equation*}
such that $[\; , \, ]_t$ is a Lie bracket. 
\end{defi}

A necessary condition for the Jacobi identity to hold is $\phi_1\in Z^2(\Lg,\Lg)$. 
The class $[\phi_1]\in H^2(\Lg,\Lg)$ is called {\it infinitesimal} deformation. 
The following result is well know.

\begin{prop}
If $H^3(\Lg,\Lg)=0$ then all obstructions vanish and
each infinitesimal deformation is integrable. 
\end{prop}

\begin{rem}
A contraction induces a formal deformation as follows.
If $\la$ contracts to $\mu$ via $g_t$, then $\la_t=g_t\cdot \la$ is
a formal deformation of $\mu$. The converse is, in general, false.
There is no duality between contractions and deformations in general.
\end{rem}

The notion of rigidity is related to algebraic deformations.  

\begin{defi}
A Lie algebra $\mu$ is called {\it formally rigid}, if every formal infinitesimal
deformation of $\mu$ is trivial. It is called {\it geometrically rigid}, if its
orbit $O(\mu)$ is open in $\CL_n(k)$. Then $\ov{O(\mu)}$ is an irreducible component
of $\CL_n(k)$. 
\end{defi}

The following result has been proved by Gerstenhaber and Schack \cite{GES}:

\begin{prop}
If $k$ has characteristic zero and $\mu$ is finite-dimensional, then
$\mu$ is geome\-trically rigid if and only if it is formally rigid.
\end{prop}

Furthermore the following results are known.

\begin{prop}
Suppose that the field is $k=\C$ or $\R$, and suppose that $H^2(\mu,\mu)=0$.
Then $\mu$ is rigid in $\CL_n(k)$.
\end{prop}

The converse is not true in general, there are explicit counter-examples.

\begin{prop}
Every complex rigid Lie algebra is algebraic. 
\end{prop}

\begin{rem}
An open question is, whether or not there is a Lie algeba law $\la \in \CN_n(k)$, which is 
rigid in the subvariety $\CN_n(k)\subset \CL_n(k)$ of nilpotent Lie algebra laws.
\end{rem}

\subsection{Generalizations}

\begin{defi}
Let $k$ be a field. A {\it discrete valuation} of $k$ is a surjective map
$\nu\colon k^{\ast}\ra \Z$ satisfying \\
\begin{align*}
\nu(xy) & = \nu(x)+\nu(y)\\
\nu(x+y) & \ge \min (\nu(x),\nu(y)).
\end{align*}
Moreover we define $\nu(0)=\infty$.
\end{defi}

The set $R=\{ x\in k \mid \nu(x)\ge 0 \} \cup \{ 0\}$ is a subring of $k$, the discrete
valuation ring (DVR) of $k$.

\begin{defi}
A {\it discrete valuation ring} is an integral domain which is the DVR of some valuation 
of its quotient field.
\end{defi}

\begin{prop}
Any discrete valuation ring is a local ring, a noetherian ring, and a pricipal ideal ring, 
hence $1$-dimensional.  
If $(t)$ is its maximal ideal, then all ideals are of the form $(t^n)$.
\end{prop}

\begin{defi}
A finitely generated extension field $K$ of $k$ of transcendence degree $1$
is called a {\it function field of dimension $1$ over $k$}.
\end{defi}

Then $K$ is a finite algebraic extension field of $k(t)$. 
Grunewald and O'Hallo\-ran proved the following theorem which shows that there
is a relationship between deformations and degenerations, see \cite{GRH}:

\begin{prop}\label{1.1.7}
Let $k$ be an algebraically closed field and $\Lg$ and $\Lg_0$ two
$n$-dimensional Lie algebras over $k$. Then $\Lg_0$ is a degeneration of $\Lg$ if and only 
if there exists a discrete valuation algebra $A$ over $k$ with quotient field $K$, and a 
Lie algebra $\La$ over $A$ of dimension $n$ such that
\begin{align}
\La \otimes_A K & \cong \Lg \otimes_k K \label{iso12} \\
\La \otimes_A k & = \Lg_0
\end{align}
\end{prop}

Note that $K$ here is a function field of dimension $1$. If $\mu_1$ represents $\Lg$ and $\mu$
represents $\La$, then \eqref{iso12} says $\phi\cdot \mu_1=\mu$ in $\CL_n(K)$, 
where $\phi\in GL(V_K)$ is the isomorphism in \eqref{iso12}. 

\begin{defi}
Let $\Lg$ be a Lie algebra over $k$ and $A$ a discrete valuation $k$-algebra with
residue field $k$. Then a Lie algebra $\La$ over $A$ is a {\it degeneration of $\Lg$ over $A$},
if there exists a finite extension $L/K$ of the quotient field $K$ of $A$, such that
\[
\La\otimes_A L \cong \Lg\otimes_k L.
\] 
The Lie algebra $\Lg_0:=\La\otimes_A k$ is called the {\it limit algebra} of the degeneration.
\end{defi}

\begin{rem}
We allow finite extensions $L/K$ in the definition. Hence we may consider also
$A$-forms of twisted versions of Lie algebras as degenerations. Note that the limit algebra
is also a degeneration in the sense of orbit closure. 
\end{rem}

\begin{defi}
Let $\mu_1\in \CL_n(k)$ and $\Lg=(V,\mu_1)$. Let $A$ be a discrete valuation $k$-algebra with 
residue field $k$ and quotient field $K$. Let $\phi \in \End(V_A)\cap GL(V_K)$.
If $\mu=\phi\cdot \mu_1$ is in $\CL_n(A)$, and hence $\mu$ defines a Lie algebra $\La=(V_A,\mu)$
over $A$, then $\La$ is called a {\it contraction} of $\Lg$ via $\phi$.
The Lie algebra $\Lg_0:=\La\otimes_A k$ is called the {\it limit algebra} of the contraction.
\end{defi}

In other words, $\Lg_0=(V,\mu_0)$ is a contraction of $\Lg=(V,\mu)$, if there is a family
$\phi_t\in \End(V_k)$ with $\det(\phi_t)\neq 0$ for $t\neq 0$, but $\det (\phi_0)=0$, such
that $\mu_0=\lim_{t\to 0}\mu_t$, where $\mu_t=\phi_t\cdot \mu_1$.

\begin{lem}
Every contraction of $\Lg$ is also a degeneration of $\Lg$, in the sense of the above
definitions.
\end{lem}

\begin{proof}
Indeed, let $\La$ be a contraction of $\Lg$ via $\phi$. Then $\La\otimes_A K$ is isomor\-phic to 
$\Lg\otimes_A K$ via $\phi$. Hence, by proposition $\ref{1.1.7}$, $\La$ is a 
degeneration of $\Lg$.
\end{proof}

\begin{prop}
A necessary condition for the existence of a contraction of $\Lg$ via $\phi$ is that 
$\phi_0(V)$ is a Lie subalgebra of $\Lg$.
In case a contraction $\Lg\ra \Lg_0$ exists, $\Lg_0$ is an extension of $\Lu=\im (\phi_0)$ 
by a nilpotent ideal $\Lv=\ker (\phi_0)$, i.e., we have a short exact sequence
\[
0 \ra \Lv \ra \Lg_0 \xrightarrow{\phi_0} \Lu \ra 0.
\]
\end{prop}

\begin{defi}
Let $\La$ be a degeneration of $\Lg$ and $\phi\colon \La\otimes_A K \ra \Lg\otimes_k K$ be 
an isomorphism. Then the pair $(\La,\phi)$ is called a {\it generalized contraction} of $\Lg$ 
with $\phi$.
\end{defi}

Hence a generalized contraction corresponds to a degeneration together with an embedding of
$\La$ in $\Lg_K$. In this sense a degeneration can be seen as a generalized contraction.

\begin{prop}
A degeneration $\La$ of $\Lg$ is isomorphic to a contraction via $\phi$ iff there exists an 
$\psi \in \Aut (\Lg_K)$ such that $\psi \kringel \phi \in \Aut(\Lg_K)\cap \End (\Lg_A)$.
\end{prop}

This means, a degeneration is a contraction, if one can choose the isomorphism in Proposition 
$\ref{1.1.7}$ from $\End(V_A)$.\\[0.5cm]
Let $A$ be a ring. We always assume that $A$ is commutative with unit.
Denote by ${\rm Spec}(A)$ the set of all proper prime ideals $\Lp$ in $A$. 
${\rm Spec}(A)$ can be turned into a topological space as follows: a subset $V$ of  ${\rm Spec}(A)$ is closed if 
and only if there exists a subset $I$ of $A$ such that $V$ consists of all those prime ideals in $A$ that contain 
$I$. This is called the Zariski topology on ${\rm Spec}(A)$. If $\Lp\in {\rm Spec}(A)$ then its residue field
is the quotient field of $A/\Lp$.

\begin{defi}
Let $\Lg_0$ be a Lie algebra over $k$ and $A$ be a $k$-algebra with specified point
$t_0\in {\rm Spec}(A)$ and residue field $k_{t_0}=k$.
A {\it deformation} of $\Lg_0$ is a Lie algebra $\La$ over $A$ together with an
isomorphism of Lie algebras over $k$,
\[
\phi\colon \Lg_0\ra \La\otimes_A k.
\]
\end{defi}

The Lie algebra $\La_k=\La\otimes_A k$ is called the {\it limit algebra} or the special fibre
of the deformation $\La$.

\begin{rem}
If $\La$ is a degeneration of $\Lg$, and $\Lg_0$ is isomorphic to the limit algebra
of $\La$ via $\phi\colon \Lg_0\ra  \La\otimes_A k$, then $(\La,\phi)$ is a deformation
of $\Lg_0$.
\end{rem}

A formal deformation of $\Lg_0$ is a deformation over the ring $A=k[[t]]$ 
of formal power series. This ring is uniquely determined as a complete 
regular $1$-dimensional local $k$-algebra.

\section{Deformations and degenerations of algebraic groups}

We want to transfer the notions to algebraic groups. Note that in the case
of Lie algebras the underlying space does not change (under degeneration or contraction). 
This will be different for algebraic groups,
where the underlying variety will also be degenerated or contracted.

\subsection{Affine group schemes}

\begin{defi}
Let $\CF$ be a sheaf of abelian groups on a topological space $X$, and $x\in X$.
Define the {\it stalk} $\CF_x$ at $x$ to be the direct limit of the abelian groups
$\CF(U)$ for all open sets $U$ containing $x$ via the restriction maps
$\rho_{UV}\colon \CF(U)\ra \CF(V)$.
\end{defi}

\begin{defi}
A {\it ringed space} is a topological space $X$ together with a sheaf of commutative rings $\CO_X$ on $X$. 
The sheaf $\CO_X$ is called the {\it structure sheaf} of $X$.
A ringed space $(X,\CO_X)$ is called a {\it locally ringed space}, if for each $x\in X$ the stalk
$\CO_{X,x}$ is a local ring. We denote by $\Lm_x$ the unique maximal ideal of $\CO_{X,x}$.
\end{defi}

\begin{defi}
Let $A$ be a ring (always commutative with unit). The {\it spectrum} of $A$ is the pair 
$({\rm Spec}(A),\CO )$ consisting of the topological space ${\rm Spec}(A)$ together 
with its structure sheaf $\CO$.
\end{defi}

If $\Lp$ is a point in  ${\rm Spec}(A)$, then the stalk $\CO_{\Lp}$ at $\Lp$ of the sheaf $\CO$ is
isomorphic to the local ring $A_{\Lp}$. Consequently,  ${\rm Spec}(A)$ is a locally ringed space.
Every sheaf of rings of this form is called an affine scheme. 

\begin{defi}
A locally ringed space $(X,\CO_X)$ is called an {\it affine scheme}, if it is
isomorphic to ${\rm Spec}(R)$ of some ring $R$, i.e., if 
\[
(X,\CO_X)\cong ({\rm Spec}(R),\CO_{{\rm Spec}(R)}).
\]
\end{defi}

\begin{ex}
If $R$ is a DVR, then ${\rm Spec}(R)$ is an affine scheme. 
\end{ex}

Its topolocial space consists of two points:
one point $t_0$ is closed, with local ring $R$. The other point $t_1$ is open and dense, with local ring
$K$, the quotient field of $R$. \\
An affine scheme $X$ is called an {\it $A$-scheme}, if its coordinate ring is an $A$-algebra.
If $Y={\rm Spec}(A)$ is an affine scheme and $\Lp\in Y$, then its residue field $k_{\Lp}$ is
the residue field of the local ring $A_{\Lp}$. 

\begin{defi}
If $X$ is an affine $A$-scheme and $\Lp\in {\rm Spec}(A)$, then the {\it fibre of $X$ over $\Lp$}
is defined by $X_{\Lp}=X\times {\rm Spec}(k_{\Lp})$. 
\end{defi}

Suppose that $A$ is a local ring with maximal ideal $\Lm$, residue field $k=k_{\Lm}$ and quotient
field $K=k_{(0)}$. Then ${\rm Spec}(A)=\{\Lm,(0) \}$.

\begin{defi}
Let $A$ be a local ring and $X$ be an $A$-scheme. 
For $\Lp=(0)$ we call the fibre $X_{\Lp}$ the {\it generic fibre of $X$} and denote it by $X_K$.
The fibre $X_{\Lm}$ is called the {\it special fibre} and is denoted by $X_k$.
\end{defi}

\begin{defi}
An {\it affine group scheme} $\CG$ over $A$ is an affine $A$-scheme $\CG$ together with 
morphisms $e\colon {\rm Spec}(A)\ra \CG$ (the identity), $i\colon \CG\ra \CG$ (the inverse), and 
$p\colon \CG\times \CG \ra \CG$ (the product), such that certain diagrams are commutative: 
Associativity, Unit and Inverse.
\end{defi}

There is the notion of a {\it smooth} affine $A$-scheme, see \cite{HA}. Note that
algebraic groups over a field $k$ of characteristic zero are smooth affine $k$-schemes. 
If we have a smooth affine group scheme $\CG$ over $A$ then we can define its Lie algebra 
$\Lie (\CG)$ via $\CG$-invariant derivations.

\subsection{Degenerations, contractions and deformations}

An affine group scheme over $A$ can be considered as a family of affine group schemes over the residue 
fields $k_t$, where $t\in {\rm Spec}(A)$. Its fibres $\CG_t$ are in fact affine group schemes with coordinate
rings $k_t[\CG_t]$. Hence we have $\CG_t=\CG_{k_t}$, and we use both notations. In particular we
write $\CG_K$ for the generic fibre of $\CG$, where $K$ is the quotient field of $A$.

\begin{defi}
Let $A$ be a discrete valuation $k$-algebra with residue field $k$ and quotient field $K$.
A {\it degeneration} of an affine algebraic group $G$ over $k$ is a smooth affine group
scheme $\CG$ over $A$, such that there is a field extension $L/K$ of finite degree, such that
$G_L$ is isomorphic to $\CG_L$.
\end{defi}

The special fiber $\CG_k$ then is called the {\it limit group} of the degeneration.

\begin{defi}
Let $A$ be a integrally closed $k$-algebra. A {\it deformation} of an affine algebraic group $G_0$ over $k$ 
is a smooth affine group scheme $\CG$ over $A$ together with a specified point  $t_0\in {\rm Spec}(A)$ and 
residue field $k_{t_0}=k$, such that there is an isomorphism of group schemes over $k$, $\psi\colon G_0 
\ra \CG_{t_0}$.
\end{defi}

\begin{defi}
Let $A$ be a discrete valuation $k$-algebra with residue field $k$ and quotient field $K$.
A {\it generalized contraction} of an affine algebraic group $G$ over $k$ 
is a pair $(\CG,\Phi)$ consisting of a degeneration $\CG$ of $G$ and an isomorphism of
$K$-group schemes $\Phi\colon \CG_K\ra G_K$. The pair $(\CG,\Phi)$ is called a {\it contraction}, 
if in addition $\Phi^{\#}(A[G])\subseteq A[\CG]$, where $\Phi^{\#}$ denotes the dual map.
\end{defi}

\begin{prop}
Let $G$ be an affine algebraic group. If $(\CG,\Phi)$ is a contraction of $G$ then 
$(\La,\phi)=(\Lie (\CG),d\Phi)$ is a contraction of $\Lg=\Lie (G)$.
\end{prop}

The same is true for a generalized contraction.

\begin{prop}
Each generalized contraction of an affine algebraic group is isomorphic to a contraction.
\end{prop}

\begin{prop}
Let $k$ be algebraically closed of characteristic zero. Then each formal degeneration of an affine algebraic group 
over $k$ (i.e., with $A=k[[t]]$) is isomorphic to a contraction.
\end{prop}

\begin{cor}
Each degeneration of a Lie algebra which corresponds to a degeneration of an affine algebraic group
is isomorphic to a contraction.
\end{cor}

\begin{defi}
Let $\La$ be a deformation or a degeneration of $\Lg$ over $A$. Then a smooth
$A$-group scheme $\CG$ with $\Lie (\CG)\cong \La$ is called a {\it lifting} of $\La$.
\end{defi}

If $G$ is an affine algebraic group over $k$ with Lie algebra $\Lg$, and if $\La$ is
a degeneration of $\Lg$ over a discrete valuation $k$-algebra $A$, then we would like
to find a lifting of $\La$ with generic fibre $G$.

\begin{defi}
Let $\La$ be a degeneration of $\Lg$ with Lie bracket $\mu=\phi\cdot \mu_1$, where
$\phi\in GL_n(K)$. A {\it conserved representation} of $\La$ is a homomorphism
\[
\rho\colon \La \ra \Lg\Ll (W_A)
\]
of $A$-Lie algebras, such that there is a $\sigma\in GL(W_K)$ with $\rho(x)=\sigma^{-1}\kringel
\rho_1(\phi(x))\kringel \sigma$ for all $x\in \Lg_K$, and such that $\rho_0=\pi_{A,k}(\rho)$
is a faithful representation of $\Lg_0=\La_k=\La\otimes_A k$.
\end{defi}

\begin{prop}[C. Daboul]
Let $\La$ be a degeneration of $\Lg$. Suppose that there  exists a conserved representation 
of $\La$, which is the derivative of a faithful representation of $G$.
Then we can construct a lifting of the degeneration.
\end{prop}

For a proof see \cite{DA}. It uses the closure of representations in the sense of schemes.
This result applies to many degenerations: if, for example, the center of the limit algebra 
is trivial, then the adjoint representation is conserved and the condition is satisfied. 
On the other hand one can use the Neron-Blowup for schemes to obtain the following
result:

\begin{prop}[C. Daboul]
All Inon\"u-Wigner contractions can be lifted to the group level.
\end{prop}

\end{document}